\documentclass[conference]{IEEEtran}
\IEEEoverridecommandlockouts
\usepackage{cite}
\usepackage{amsmath, amssymb, amsfonts, amsthm, mathtools}
\usepackage{algorithmic, algorithm}
\usepackage{graphicx}
\usepackage{textcomp}
\usepackage{xcolor}
\usepackage{hyperref}

\usepackage{xspace} 

\newcommand{\algname}{\textsc{DT-GO}\xspace} 

\newcommand{\norm}[1]{\left\lVert#1\right\rVert}
\newcommand{\normsq}[1]{\left\lVert#1\right\rVert^2}
\newcommand{\expec}[1]{\mathbb{E}\left[#1\right]}

\newcommand{\R}{\mathbb{R}}

\DeclareMathOperator{\diag}{diag}

\usepackage{enumerate, cleveref}
\crefname{assumption}{Assumption}{Assumptions}
\crefname{figure}{Fig.}{Figs.}
\crefname{definition}{Def.}{Defs.}
\crefname{equation}{Eq.}{Eqs.}
\crefname{algorithm}{Alg.}{Algs.}
\crefname{lemma}{Lemma}{Lemmas}
\crefname{theorem}{Theorem}{Theorems}
\crefname{section}{Section}{Sections}
\crefname{proposition}{Proposition}{Propositions}

\theoremstyle{plain}
\newtheorem{theorem}{Theorem}[section]
\newtheorem{proposition}[theorem]{Proposition}
\newtheorem{lemma}[theorem]{Lemma}

\theoremstyle{definition}
\newtheorem{definition}{Definition}
\newtheorem{assumption}[theorem]{Assumption}

\theoremstyle{remark}

\usepackage{multirow, booktabs} 

\usepackage{pgfplots}
  \pgfplotsset{compat=newest}
  \usetikzlibrary{plotmarks}
  \usetikzlibrary{arrows.meta}
  \usetikzlibrary{positioning}
  \usepgfplotslibrary{patchplots}
  \usepackage{grffile}

\begin{document}

\title{Decentralized Optimization in Networks with Arbitrary Delays
    \thanks{This work was supported in part by the NSF Award ECCS-2207457.}
}

\author{\IEEEauthorblockN{Tomàs Ortega, Hamid Jafarkhani}
    \IEEEauthorblockA{\textit{Center for Pervasive Communications and Computing} \\
        \textit{University of California, Irvine}\\
        Irvine, USA \\
        Emails: tomaso@uci.edu, hamidj@uci.edu}
}

\maketitle

\begin{abstract}
    We consider the problem of decentralized optimization in networks with communication delays.
    To accommodate delays, we need decentralized optimization algorithms that work on directed graphs.
    Existing approaches require nodes to know their out-degree to achieve convergence.
    We propose a novel gossip-based algorithm that circumvents this requirement, allowing decentralized optimization in networks with communication delays.
    We prove that our algorithm converges on non-convex objectives, with the same main complexity order term as centralized Stochastic Gradient Descent (SGD), and show that the graph topology and the delays only affect the higher order terms.
    We provide numerical simulations that illustrate our theoretical results.
\end{abstract}


\section{Introduction}
\begin{figure}[htbp]
    \centering
    \tikzstyle{edge} = [very thick, -{Latex}]
\tikzstyle{vedge} = [very thick, magenta, dashdotted, -{Latex}]
\tikzstyle{non-virtualnode} = [draw, circle, fill=white, very thick]
\tikzstyle{virtualnode} = [draw, circle, magenta, fill=magenta!5, very thick, dashdotted]

\begin{tikzpicture}

    \node[non-virtualnode] (3) {\textbf{3}};
    \node[non-virtualnode] (4) [left = 1cm of 3] {\textbf{4}};
    \node[non-virtualnode] (1) [below = 1cm of 3] {\textbf{1}};
    \node[non-virtualnode] (2) [below = 1cm of 4] {\textbf{2}};
    \node[non-virtualnode] (5) [right = 1cm of 3] {\textbf{5}};

    \draw[edge] (3) -- (4);
    \draw[edge] (4) -- (2);
    \draw[edge] (2) -- (3);
    \draw[edge] (2) -- (1);
    \draw[edge] (1) -- (3);
    \draw[edge] (5) to[bend left=10] (3);
    \draw[edge] (3) to[bend left=10] (5);

    \draw[edge] (3) to[loop above, looseness=10] (3);
    \draw[edge] (4) to[loop left, looseness=10] (4);
    \draw[edge] (1) to[loop right, looseness=10] (1);
    \draw[edge] (2) to[loop left, looseness=10] (2);
    \draw[edge] (5) to[loop right, looseness=10] (5);

\end{tikzpicture}
    \caption{Directed communication graph example, with $N=5$.}
    \label{fig:example-digraph}
\end{figure}
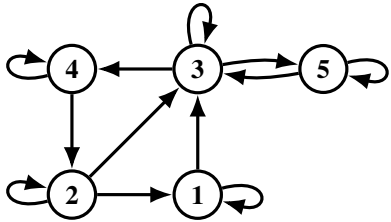

In decentralized optimization, the nodes in a network cooperate to minimize a global objective function that is the average of nodes' local objective functions.
The motivation for these problems comes from a variety of applications, including, but not limited to, decentralized estimation in sensor networks, collaborative machine learning, and decentralized coordination of multi-agent systems.
We formalize these problems by considering a collection of $N$ nodes, which can communicate with each other through the edges of a graph $G$ -- see \cref{fig:example-digraph} for an example.
Each node has a local cost function $f_n : \R^d \to \R$.
The goal is for all nodes to find a common model $x \in \R^d$ that minimizes a global cost function $f$ of the form
\begin{equation}  \label{eq:objective-function}
    f(x) \coloneqq \frac{1}{N} \sum_{n=1}^N f_n(x).
\end{equation}
Previous research has proposed several techniques to minimize \eqref{eq:objective-function}.
The first study of this problem can be tracked back at least to \cite{tsitsiklis1984problems}.
For the problem of finding a common solution among nodes, the seminal works in gossip algorithms \cite{Kempe_Dobra_Gehrke_2003,gossip_algorithms} brought ideas from mixing in Markov chains to allow averaging over graphs.
To minimize \eqref{eq:objective-function}, it is standard to either use a combination of gradient descent with gossip steps \cite{distributed_subgradient}, or problem-specific methods such as alternating direction method of multipliers (ADMM) \cite{Wei_Ozdaglar_2012}.
In the case where $G$ is undirected, there are multiple efforts to solve this problem in various settings, like asynchronous communications \cite{Srivastava_Nedic_2011}, time-varying graphs \cite{time-varying-undirected-graphs, unified_decentralized_SGD}, or quantized communications \cite{chocosgd, ortegaGossip}.
In the case where $G$ is directed, nodes must know their out-degree to achieve consensus, as stated in \cite{decentralized_optimization_survey_nedic}.
Under such assumption, a family of so-called Push-Sum methods was proposed, starting with \cite{weighted_gossip,distributed_subgradient}.
There are also multiple extensions for various settings, like quantized communications \cite{chocosgd}, and asynchronous communications \cite{SGP-distributed}.

There exist several nice surveys that detail the contributions in this field, see for example \cite{decentralized_optimization_survey_nedic} for synchronous methods, and \cite{advances_in_asynchronous_optimization} for asynchronous methods.

The common assumption for existing algorithms for directed decentralized optimization is that nodes in the network must know their out-degree.
In this paper, we propose a gossip-based optimization approach where nodes \emph{do not} need to know their out-degree to achieve convergence.
Such a property is desired in many scenarios which arise naturally.
For example, consider networks where links are directed and established without a handshake, in other words, networks where the transmitter sends messages and the receivers do not need to acknowledge their reception.
Such a setting is likely in scenarios where some nodes may have a transmit power that is larger than others, making some bidirectional communication impossible.
Another possible application is in networks where nodes simply broadcast messages periodically, without knowledge of who receives them.
Our algorithm, \algname (Delay Tolerant Gossiped Optimization), allows decentralized optimization in such networks.

\section{Setup and proposed algorithm}
\subsection{Problem setup}
As described in the introduction, we consider a set of $N$ nodes that share their information over a network modeled as a \emph{directed} graph $G$.
Since the shared information can only reach a limited number of neighboring nodes, not necessarily known to the transmitter, we assume that the nodes can \emph{only} know their \emph{in-degree}.

Our approach to minimize \eqref{eq:objective-function} is to iteratively perform two phases: (i) local optimization, and (ii) consensus.
At each local optimization phase, nodes optimize their models based on their local data.
This makes local models drift from the average.
Therefore, we introduce a consensus phase, where nodes communicate among themselves to converge to the average of their models.
The problem in the consensus phase is known in the literature as decentralized averaging.

\subsection{Decentralized averaging} \label{sec:decentralized-averaging}
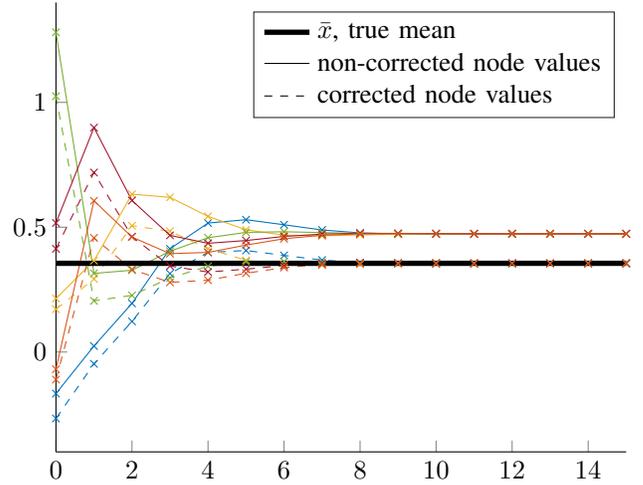
\begin{figure}[htbp]
    \centering
%
%
\definecolor{mycolor1}{rgb}{0.00000,0.44700,0.74100}%
\definecolor{mycolor2}{rgb}{0.92900,0.69400,0.12500}%
\definecolor{mycolor3}{rgb}{0.46600,0.67400,0.18800}%
\definecolor{mycolor4}{rgb}{0.63500,0.07800,0.18400}%
\definecolor{mycolor5}{rgb}{0.85000,0.32500,0.09800}%
\begin{tikzpicture}

\begin{axis}[%
width=0.951 * 0.9\columnwidth,
height=0.75 * 0.9\columnwidth,
at={(0\columnwidth,0\columnwidth)},
scale only axis,
xmin=0,
xmax=15,
ymin=-0.4,
ymax=1.4,
axis background/.style={fill=white},
axis x line*=bottom,
axis y line*=left,
legend style={legend cell align=left, align=left, draw=white!15!black}
]
\addplot [color=mycolor1, mark=x, mark options={solid, mycolor1}, forget plot]
  table[row sep=crcr]{%
0	-0.166382762898544\\
1	0.0240198738971338\\
2	0.194856901964998\\
3	0.413578696792329\\
4	0.516599869140531\\
5	0.530029841317183\\
6	0.50976929817379\\
7	0.488927533177665\\
8	0.477216912919211\\
9	0.473072172825798\\
10	0.47260516474605\\
11	0.473148039672441\\
12	0.473600420868596\\
13	0.473768224524826\\
14	0.473763100213317\\
15	0.473708365188305\\
};
\addplot [color=mycolor1, dashed, mark=x, mark options={solid, mycolor1}, forget plot]
  table[row sep=crcr]{%
0	-0.266212420637672\\
1	-0.0473372060417108\\
2	0.122608968992289\\
3	0.314224681144008\\
4	0.39907658842509\\
5	0.406624924310868\\
6	0.386979866136379\\
7	0.368312321899792\\
8	0.358199594948836\\
9	0.354770269136197\\
10	0.354470623185172\\
11	0.354975303510421\\
12	0.355364766909793\\
13	0.355498715211849\\
14	0.355485821005468\\
15	0.355434668483628\\
};
\addplot [color=mycolor2, mark=x, mark options={solid, mycolor2}, forget plot]
  table[row sep=crcr]{%
0	0.214422510692812\\
1	0.365693930032861\\
2	0.632300491619661\\
3	0.619621041488732\\
4	0.543459813493836\\
5	0.489508755030396\\
6	0.46808576818154\\
7	0.465506292660758\\
8	0.468927432732385\\
9	0.472138156666302\\
10	0.473690914598832\\
11	0.47405280206475\\
12	0.473936028181057\\
13	0.473757975901807\\
14	0.473653630163293\\
15	0.47361988392922\\
};
\addplot [color=mycolor2, dashed, mark=x, mark options={solid, mycolor2}, forget plot]
  table[row sep=crcr]{%
0	0.17153800855425\\
1	0.292555144026288\\
2	0.505840393295726\\
3	0.483928495706172\\
4	0.414173260196646\\
5	0.367334807961891\\
6	0.349644777663205\\
7	0.34808686799788\\
8	0.351340943323558\\
9	0.354170977234147\\
10	0.35547998383567\\
11	0.355754230309165\\
12	0.355632663513905\\
13	0.355472926799086\\
14	0.355383515961788\\
15	0.355356392012518\\
};
\addplot [color=mycolor3, mark=x, mark options={solid, mycolor3}, forget plot]
  table[row sep=crcr]{%
0	1.28084875704001\\
1	0.314976129509147\\
2	0.327655579640076\\
3	0.403816807634972\\
4	0.457767866098412\\
5	0.479190852947268\\
6	0.48177032846805\\
7	0.478349188396422\\
8	0.475138464462506\\
9	0.473585706529976\\
10	0.473223819064057\\
11	0.473340592947751\\
12	0.473518645227001\\
13	0.473622990965515\\
14	0.473656737199588\\
15	0.473655960440514\\
};
\addplot [color=mycolor3, dashed, mark=x, mark options={solid, mycolor3}, forget plot]
  table[row sep=crcr]{%
0	1.024679005632\\
1	0.204907553668069\\
2	0.226819451257624\\
3	0.29657468676715\\
4	0.343413139001905\\
5	0.36110316930059\\
6	0.362661078965915\\
7	0.359407003640238\\
8	0.356576969729649\\
9	0.355267963128126\\
10	0.354993716654631\\
11	0.355115283449891\\
12	0.355275020164709\\
13	0.355364431002008\\
14	0.355391554951278\\
15	0.355389545203905\\
};
\addplot [color=mycolor4, mark=x, mark options={solid, mycolor4}, forget plot]
  table[row sep=crcr]{%
0	0.516965349372911\\
1	0.898907053206461\\
2	0.606941591357804\\
3	0.46729858549894\\
4	0.435557696566956\\
5	0.446662781332684\\
6	0.462926817139976\\
7	0.472348572804013\\
8	0.475348880600218\\
9	0.475243672531362\\
10	0.474414689530669\\
11	0.473819254297363\\
12	0.473579923622557\\
13	0.473549284424779\\
14	0.473586137695147\\
15	0.473621437447367\\
};
\addplot [color=mycolor4, dashed, mark=x, mark options={solid, mycolor4}, forget plot]
  table[row sep=crcr]{%
0	0.413572279498326\\
1	0.719125642565165\\
2	0.462016598116617\\
3	0.344418024687121\\
4	0.320496355727135\\
5	0.33195474736452\\
6	0.346528958332555\\
7	0.354595018649235\\
8	0.357001011144737\\
9	0.356788990437193\\
10	0.356028476782659\\
11	0.355511096718645\\
12	0.355313190084268\\
13	0.355294105124489\\
14	0.355329268063248\\
15	0.355360411507263\\
};
\addplot [color=mycolor5, mark=x, mark options={solid, mycolor5}, forget plot]
  table[row sep=crcr]{%
0	-0.0689839867976905\\
1	0.60593238512116\\
2	0.460454257315153\\
3	0.394054918477614\\
4	0.398935863056293\\
5	0.428351864577353\\
6	0.45377135876231\\
7	0.46777084361518\\
8	0.473060016005801\\
9	0.474099240234154\\
10	0.473842473382065\\
11	0.473533146223061\\
12	0.473436869585406\\
13	0.473477757406203\\
14	0.473550374185859\\
15	0.473603555692724\\
};
\addplot [color=mycolor5, dashed, mark=x, mark options={solid, mycolor5}, forget plot]
  table[row sep=crcr]{%
0	-0.110374378876306\\
1	0.457152313377849\\
2	0.331029933522959\\
3	0.278924692390292\\
4	0.287749689578721\\
5	0.315581414290313\\
6	0.338342291795452\\
7	0.350501685380684\\
8	0.354954344510461\\
9	0.355765657120055\\
10	0.35551681012409\\
11	0.355255263389361\\
12	0.355185273419626\\
13	0.355230146792168\\
14	0.355297288897088\\
15	0.355344421924183\\
};
\addplot [color=black, line width=2.0pt]
  table[row sep=crcr]{%
0	0.3553739734819\\
1	0.3553739734819\\
2	0.3553739734819\\
3	0.3553739734819\\
4	0.3553739734819\\
5	0.3553739734819\\
6	0.3553739734819\\
7	0.3553739734819\\
8	0.3553739734819\\
9	0.3553739734819\\
10	0.3553739734819\\
11	0.3553739734819\\
12	0.3553739734819\\
13	0.3553739734819\\
14	0.3553739734819\\
15	0.3553739734819\\
};
\addlegendentry{$\bar x$, true mean}

\addplot [color=black]
  table[row sep=crcr]{%
0	0.3553739734819\\
};
\addlegendentry{non-corrected node values}

\addplot [color=black, dashed]
  table[row sep=crcr]{%
0	0.3553739734819\\
};
\addlegendentry{corrected node values}

\end{axis}
\end{tikzpicture}%
    \caption{Plot of corrected and non-corrected node values throughout time. Initial node values are chosen at random from a normal $\mathcal{N}(0,5)$. The gossip weights are the inverse of the in-degrees, where $G$ is shown in \cref{fig:example-digraph}.}
    \label{fig:corrected-values}
\end{figure}

Let us assume that every Node $n$ starts with an initial model, or initial state, $x_n(0) \in \R^d$.
At each time $t$, nodes broadcast their models and collect the models of their neighbors.
They then perform a weighted average with the models they receive, which is called a \emph{gossip} step.
Specifically, Node $n$ weighs information received from Node $m$ with weight $W_{nm}$.
A natural choice for these weights is the inverse of the in-degree of each node.
Other choices also work, as long as they satisfy our definition of a gossip matrix.
\begin{definition}[Gossip matrix] \label{def:gossip-matrix}
    Given a strongly connected directed graph $G$, its gossip matrix $W$ is an $N\times N$ real matrix whose entries satisfy:
    \begin{enumerate}[(i)]
        \item $W$ is row-stochastic: the sum of the entries in every row is one, i.e., $\sum_{j=1}^N W_{ij} = 1$ for all $i$.
        \item Every entry $W_{nm}$ is non-negative, and zero only if there is no directed edge from Node $m$ to Node $n$ in $G$.
        \item Entries in the diagonal are positive. This is equivalent to asking each node in $G$ to have a self-loop.
    \end{enumerate}
\end{definition}
Formally, each gossip iteration is defined as $x_n(t+1) = \sum_{m=1}^N W_{nm} x_m(t)$, which can be written in a matrix form as $X(t+1) = W X(t)$.
Here, $X(t)$ is the $N \times d$ real matrix whose rows are the node states.
Notice from item (ii) in \cref{def:gossip-matrix} that $W_{nm}$ is only positive when Node $n$ can actually receive model $x_n(t)$.
When $t \to \infty$, we want the sequence of models to tend to the average, i.e, $x_n(t) \to \bar x = \frac{1}{N} \sum_{n=1}^N x_n(0)$.
However, if the gossip matrix $W$ is row-stochastic nodes \emph{do not} necessarily converge to the average. 
An illustration of this phenomenon can be seen in \cref{fig:corrected-values}, where the nodes using a regular gossip algorithm, marked by the ``non-corrected node values'', clearly do not converge to $\bar x$.
For row-stochastic matrices, it can be shown that the nodes converge to a weighted average
\begin{equation} \label{eq:weighted-average}
    \tilde x = \sum_{n=1}^N \pi_n x_n(0),
\end{equation}
where $\pi_n$ are non-negative weights that add up to one \cite[Lemma 5]{nedic_consensus}.
We prove that $\pi_1, \ldots, \pi_N$ are positive in \cref{sec:analysis}.
In this section, first, we present an algorithm that ``corrects'' the values such that the nodes converge to the average.
Then, we show that this algorithm is capable of handling links with communication delays.

\subsection{Algorithm design} \label{sec:algorithm-design}
We present \algname, which corrects the weighted average described in \cref{eq:weighted-average} and allows the consensus phase to converge to the true average.
Then, we can minimize the cost function in \cref{eq:objective-function} judiciously by interchanging consensus and local optimization phases. 

The key idea is described in the following Lemma.
\begin{lemma}
    Consider a digraph $G$ with an associated gossip matrix $W$ as in \cref{def:gossip-matrix}.
    If every node $n$ multiplies its initial state $x_n(0)$ by a factor $d_n = \frac{1}{N \pi_n}$, then the gossip iterations $x_n(t+1) = \sum_{m=1}^N W_{nm} x_m(t)$ converge to the true mean, i.e.,
    \begin{equation}
        \sum_{n=1}^N \pi_n d_n x_n(0) = \sum_{n=1}^N \frac{1}{N} x_n(0) = \bar x.
    \end{equation}
\end{lemma}
The proof of this lemma follows immediately from \cref{eq:weighted-average} by inserting the re-weighted initial states.
See \cref{fig:corrected-values} for an illustration of this lemma in practice, where the corrected node values correspond to a run where every initial state from Node $n$ is re-weighted with $d_n$.

Observe that if every node has initial state $x_n(0) = e_n$, a real vector of length $N$ that includes all zeros, except a one in coordinate $n$, the gossip iterations converge to
\begin{equation} \label{eq:obtaining-pin}
    \tilde x = \sum_{n=1}^N \pi_n e_n,
\end{equation}
that is, a vector with each $\pi_n$ in the $n$-th coordinate.
Therefore, each node can compute $d_n$ when convergence is achieved.
Leveraging this observation, we propose \cref{alg:decentralized_opt_alg}.
\begin{algorithm}[b!] 
    \caption{\algname at Node $n$}\label{alg:decentralized_opt_alg}
    \begin{algorithmic}[1]
        \STATE Generate local id number $id_n$. \\ \COMMENT {\textsc{Warm-Up Period}}
        \STATE {Initialize dictionary $dict \gets \{id_n:1\}$.}
        \FOR{$T_{\text{warm-up}}$ rounds}
        \STATE {Broadcast $dict$ and receive neighbors' dictionaries.}
        \STATE {$dict \gets$ weighted average of available dictionaries.}
        \ENDFOR{}
        \STATE{From $dict$, obtain $N$ and $\pi_n$, as in \cref{eq:obtaining-pin}. \\ \COMMENT {\textsc{Minimizing $f$ Period}}}
        \STATE {Initialize $x_n(0)$.}
        \FOR{$t$ in $0, \ldots, T-1$}
        \STATE {Initialize auxiliary variables $z_n, y$. \\ \COMMENT {Optimization phase}}
        \STATE {SGD step: $y \gets x_n(t) - \eta \nabla F_n(x_n(t), \xi_n)$. \\ \COMMENT {Consensus phase}}
        \STATE {Adjust update: $z_n \gets x_n(t) + \frac{1}{N\pi_n} (y - x_n(t))$.}
        \FOR {$\tau_g$ rounds}
        \STATE {Broadcast $z_n$ and receive neighbor states.}
        \STATE {Average states: $z_n \gets \sum_{m=1}^N W_{nm} z_m$.}
        \ENDFOR{}
        \STATE {$x_n(t+1) \gets z_n$.}
        \ENDFOR{}
    \end{algorithmic}
\end{algorithm}
Note that we do not use vectors $e_n$, because nodes do not know the network size $N$ a priori.
Instead, we use dictionaries that follow the same principle, and are implementable in practice.
Once the warm-up period is over, nodes can obtain $N$ from the dictionary size.

\subsection{Incorporating Delays} \label{sec:incorporating-delays}
The framework proposed for \algname is designed to easily accommodate links with delays, i.e., links where information takes more than one round to arrive at the receiver.
To incorporate the delays, similar to \cite{nedic_consensus, tsitsiklis1984problems}, we introduce the notion of virtual nodes, or non-computing nodes, that serve as a relay of the message for a round.
If Node $n$ sends messages to Node $m$ with a delay of $k$ rounds, we simply add $k$ nodes to the graph $G$, each with objective function 0 to avoid modifying the objective function \eqref{eq:objective-function}.
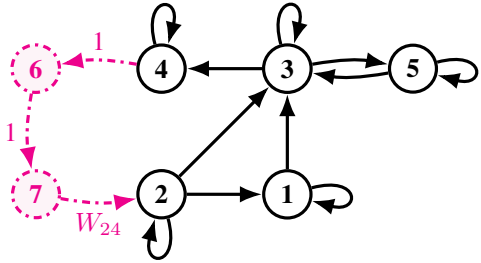
\begin{figure}[htbp]
    \centering
    \tikzstyle{edge} = [very thick, -{Latex}]
\tikzstyle{vedge} = [very thick, magenta, dashdotted, -{Latex}]
\tikzstyle{non-virtualnode} = [draw, circle, fill=white, very thick]
\tikzstyle{virtualnode} = [draw, circle, magenta, fill=magenta!5, very thick, dashdotted]

\begin{tikzpicture}

    \node[non-virtualnode] (3) {\textbf{3}};
    \node[non-virtualnode] (4) [left = 1cm of 3] {\textbf{4}};
    \node[non-virtualnode] (1) [below = 1cm of 3] {\textbf{1}};
    \node[non-virtualnode] (2) [below = 1cm of 4] {\textbf{2}};
    \node[non-virtualnode] (5) [right = 1cm of 3] {\textbf{5}};
    \node[virtualnode] (6) [left = 1cm of 4] {\textbf{6}};
    \node[virtualnode] (7) [left = 1cm of 2] {\textbf{7}};

    \draw[vedge] (4) to[bend right=10] node[above]{1} (6);
    \draw[vedge] (6) to[bend right=10] node[left]{1} (7);
    \draw[vedge] (7) to[bend right=10] node[below]{$W_{24}$} (2);
    \draw[edge] (3) -- (4);
    \draw[edge] (2) -- (3);
    \draw[edge] (2) -- (1);
    \draw[edge] (1) -- (3);
    \draw[edge] (5) to[bend left=10] (3);
    \draw[edge] (3) to[bend left=10] (5);

    \draw[edge] (3) to[loop above, looseness=10] (3);
    \draw[edge] (4) to[loop above, looseness=10] (4);
    \draw[edge] (1) to[loop right, looseness=10] (1);
    \draw[edge] (2) to[loop below, looseness=10] (2);
    \draw[edge] (5) to[loop right, looseness=10] (5);

\end{tikzpicture}
    \caption{Example of a graph with delayed links. We have added a delay of 2 rounds to the edge between Node 4 and Node 2 of the graph in \cref{fig:example-digraph}.}
    \label{fig:delayed-graph}
\end{figure}
For example, in \cref{fig:delayed-graph}, Node 4 sends messages to Node 2 with a delay of 2 rounds.
Therefore, we add two nodes, 6 and 7, to the graph, and Node 4 sends messages to Node 6, which sends messages to Node 7, which finally sends messages to Node 2.
The weight of the old edge from Node 4 to 2, $W_{24}$, is now assigned to the edge from Node 7 to 2 and the weights of the other new edges are set to 1.
Note that this procedure ensures that the nodes' input weights still add up to one.
Also, the resulting extended graph is still strongly connected.

\section{Performance bounds and proof of convergence} \label{sec:analysis}
In this section, first, we discuss the case without delays, as it gives a stronger and more illustrative bound.
Then, we present the case with delays.
\subsection{Case without delays}
First, we introduce a proposition that is needed to show that our algorithm converges.
\begin{proposition} \label{prop:gossip-properties}
    Given a gossip matrix $W$ satisfying \cref{def:gossip-matrix} we can ensure that
    \begin{enumerate}[(i)]
        \item The limit $\lim_{t \to \infty} W^t \coloneqq W^\infty$ exists, that is, our gossip algorithm converges to a stationary solution.
        \item The matrix $W^\infty$ is row-stochastic, and its rows are all identical, with positive entries $\pi_1, \ldots, \pi_N$ that add up to one.
        \item The squared Frobenius norm of the global state converges to the stationary solution at a geometric rate, i.e., there exist constants $C$ and $\rho$, with $\rho < 1$, such that for all $t>0$, $\norm{W^t - W^\infty}_2^2 \leq C \rho^t$.
    \end{enumerate}
\end{proposition}
\begin{IEEEproof}
    Since $W$ is primitive, statements (i) and (ii) are immediate consequences of applying the Perron-Frobenius theorem \cite{Meyer_Stewart_2023}.
    Furthermore, $\pi_1, \ldots, \pi_N$ are the coordinates of the left Perron eigenvector of $W$. 
    Statement (iii) also has a straightforward proof via eigendecomposition, as shown in \cite[Fact 3]{convergence_rate_markov}.
    In addition, we have $\rho = |\lambda_2|^2$, the second largest eigenvalue of $W$ in absolute value squared.
\end{IEEEproof}

\subsection{Case with delays}
Let us now assume that we have a network with arbitrary delays.
The number of non-virtual nodes is $N$, and we add virtual nodes to obtain a network of size $N_v$.
Similarly, the original gossip matrix $W$ is extended as previously described to a gossip matrix $W_v$.
\begin{proposition} \label{prop:gossip-properties-with-delays}
    Given a gossip matrix $W$ satisfying \cref{def:gossip-matrix}, and its extended version with delays $W_v$, we can ensure that
    \begin{enumerate}[(i)]
        \item The limit $\lim_{t \to \infty} W_v^t \coloneqq W_v^\infty$ exists, that is, our gossip algorithm with delays converges to a stationary solution.
        \item The matrix $W_v^\infty$ is row-stochastic, and its rows are all identical, with non-negative entries $\pi_{v,1}, \ldots, \pi_{v,N_v}$ that add up to one.
        \item The squared Frobenius norm of the global state converges to the stationary solution at a geometric rate, i.e., there exist constants $C$ and $\rho$, with $\rho < 1$, such that for all $t>0$, $\norm{W_v^t - W_v^\infty}_2^2 \leq C \rho^t$.
        \item The weights that correspond to non-virtual nodes, i.e., $\pi_{v,1}, \ldots, \pi_{v,N}$ are all positive.
    \end{enumerate}
\end{proposition}
\begin{IEEEproof}
    From \cite[Lemma 5]{nedic_consensus}, Facts (i) and (ii) follow immediately.
    For Fact (iii), \cite[Lemma 5]{nedic_consensus} ensures that
    \begin{equation}
        \norm{W^t - W^\infty}_F \leq 2 \frac{1 + \eta^{-B_2}}{1 - \eta^{B_2}}(1 - \eta^{B_2})^{\frac{t}{B_2}},
    \end{equation}
    where $B_2 \coloneqq N - 1 + NB_1$, $B_1$ is the maximum number of delays between two nodes in $G$, and $\eta \in (0,1)$ is a positive lower bound for all the non-zero entries of $W$.
    Squaring the expression and knowing that the Frobenius norm upper-bounds the 2-norm proves the statement.
    The proof of (iv) follows from \cite[Lemma 5.2.1]{tsitsiklis1984problems}.
\end{IEEEproof}
\cref{prop:gossip-properties-with-delays} allows us to guarantee convergence for graphs with delays, albeit it often gives us a looser bound for delay-less graphs.

\subsection{\algname's convergence rate}
A corollary of \cref{prop:gossip-properties,prop:gossip-properties-with-delays} is that we can choose a number of gossip rounds $\tau_g$ such that $C \rho^{\tau_g} < 1$, as needed for our convergence results.

Note that if $W$ is doubly stochastic and we do not add delays, $C=1$, and $\tau_g$ can be set to 1.
Also, $W$ then converges to an all-one matrix divided by $N$, and our algorithm becomes equivalent to Decentralized SGD (DSGD) \cite{distributed_subgradient}.
Furthermore, if $G$ is a complete graph and we set edge weights to $1/N$, then the gossip matrix becomes the averaging matrix, and both \algname and DSGD are equivalent to centralized SGD.

Assuming that the warm-up phase is complete and $\tau_g$ is large enough such that $C \rho^{\tau_g} < 1$, with a slight abuse of the notation we can write \cref{alg:decentralized_opt_alg}'s iterations as
\begin{equation} \label{eq:alg-iterations}
    X(t+1) = W^{\tau_g}\left(X(t) - \eta D \partial F(t) \right),
\end{equation}
where $\partial F(t)$ is an $N \times d$ real matrix that represents the stochastic gradients at time $t$, and $D = \diag(\tfrac{1}{N\pi_1}, \ldots, \tfrac{1}{N\pi_N})$ is the diagonal correction matrix.
If there are virtual nodes, the corresponding entries in the diagonal matrix $D$ are set to $1$.
This ensures that $D$ is well defined, since $\pi_n > 0$ when Node $n$ is non-virtual -- see \cref{prop:gossip-properties-with-delays}.
This matrix notation allows us to follow an analysis approach similar to \cite{decentralized_DL_with_compression}.
Before stating our main convergence result, first, we make the following standard assumptions \cite{chocosgd,distributed_subgradient}.
\begin{assumption}\label{ass:f}
    Each function $f_n \colon \R^d \to \R$ is $L$-smooth, that is
    \begin{equation}
        \norm{\nabla f_n(y)-\nabla f_n(x)} \leq L \norm{y-x}\, \forall x,y \in \R^d,
    \end{equation}
    and we have access to unbiased stochastic gradients, $\nabla F_n(x, \xi_n)$, with bounded variance on each worker:
    \begin{align} \label{eq:assump_second_momentum}
        \expec{\normsq{\nabla F_n(x, \xi_n) - \nabla f_n(x)}} & \leq \sigma_n^2\,                       \\
        \expec{\normsq{\nabla F_n(x, \xi_n)}}                 & \leq G^2\,        & \forall x \in \R^d.
    \end{align}
    We denote $\overline{\sigma}^2 := \frac{1}{N}\sum_{n = 1}^N\sigma_n^2$ for convenience.
\end{assumption}
We now present our main convergence result, which applies for cases with and without delays.
\begin{theorem}\label{th:non-convex-sigma}
    Given a gossip matrix as in \cref{def:gossip-matrix} (possibly extended with delays as in \cref{sec:incorporating-delays}), a number of iterations $T$, a large enough $\tau_g$ such that $\norm{W^{\tau_g} - W^\infty}^2_2 \leq C\rho^{\tau_g} \coloneqq 1-c < 1$, and under \cref{ass:f}, there exists a constant (depending on $T$) stepsize $\eta$ such that the \algname iterations $\tilde{x}^{(t)} \coloneqq \sum_{n=1}^N \pi_n x_n(t)$ satisfy:
    \begin{multline*}
        \frac{1}{T}\sum_{t = 0}^{T-1} \expec{\normsq{\nabla f \left(\tilde{x}(t) \right)}_2} = \\
        {} \mathcal{O} \left( \left( \frac{L F_0 \overline{\sigma}^2 }{N T} \right)^{1/2} +  \left(\frac{ \norm{D}_2 G L F_0}{c T}\right)^{2/3} + \frac{L F_0}{T}\right),
    \end{multline*}
    where $F_0 := f(\tilde{x}^{(0)}) - f^\star$, $f^\star$ is a lower bound for $f$, and $D$ is the diagonal correction matrix, with virtual nodes having 1 in their corresponding entries.
\end{theorem}
We present a sketch of \cref{th:non-convex-sigma}'s proof, which is not included with full details due to the lack of space.
The first key insight is that
\begin{equation} \label{eq:insight}
    \expec{f(\tilde x (t+1)} = \expec{f(\tilde x(t)) - \frac{\eta}{N} \sum_{n=1}^N \nabla F_n(x_n(t), \xi_n)},
\end{equation}
since the weighted averages are preserved at each iteration.
We can check this as follows: let us define $\tilde X(t+1)$ as an $N \times d$ matrix whose rows are $\tilde x(t)$.
From \cref{eq:alg-iterations}, we obtain
\begin{align*}
    \tilde X(t+1) = W^\infty X(t+1) & = W^\infty \left(X(t) - \eta D \partial F(t) \right) \\
                                    & = \tilde X(t) - \eta J \partial F(t),
\end{align*}
where $J$ is the averaging matrix of all ones divided by $N$.
Notice that this is exactly \cref{eq:insight} in matrix form.
Using this insight, we expand \cref{eq:insight} leveraging $L$-smoothness:
\begin{align*}
    \expec{f(\tilde x (t+1))} & \leq \expec{f(\tilde x(t))}                                                                              \\
                              & - \expec{ \langle \nabla f (\tilde x(t)), \frac{\eta}{N} \sum_{n=1}^N \nabla F_n(x_n(t), \xi_n) \rangle} \\
                              & + \eta^2 \frac{L}{2}\expec{ \normsq{\frac{1}{N} \sum_{n=1}^N \nabla F_n(x_n(t), \xi_n) }_2}.
\end{align*}
After algebraic manipulation, and using the bounded variance assumption as well as $L$-smoothness, we obtain
\begin{align*}
    \expec{f(\tilde x (t+1))} & \leq \expec{f(\tilde x(t))} - \frac{\eta}{4}\expec{ \normsq{\nabla f(\tilde x(t) )}_2 }                              \\
                              & +\frac{3\eta L^2}{4N}\sum_{n=1}^N \expec{ \normsq{\tilde x(t) - x_n(t)}_2 }+ \eta^2\frac{L \overline{\sigma}^2}{2N}.
\end{align*}
Rearranging terms and averaging over $T$ iterations yield
\begin{align*}
    \frac{1}{T} \sum_{t=0}^{T-1} \expec{ \normsq{\nabla f(\tilde x(t) )}_2 } & \leq 4\frac{f(\tilde x(0)) - f^\star}{\eta T}  + \eta\frac{2L \overline{\sigma}^2}{N}     \\
                                                                             & +\frac{3L^2}{NT} \sum_{t=0}^{T-1} \sum_{n=1}^N \expec{ \normsq{\tilde x(t) - x_n(t)}_2 }.
\end{align*}
We bound the third term on the RHS using \cref{lem:bound_general}.
Then, choosing the appropriate stepsize, following arguments similar to \cite[Lemma A.5]{decentralized_DL_with_compression} concludes the proof.

Notice that the first term in \cref{th:non-convex-sigma} is the same as that of centralized SGD, and only the higher order terms are affected by the graph topology and the delays between nodes, which affect $D$ and $c$ as described in \cref{prop:gossip-properties-with-delays}.

We now present \cref{lem:bound_general}, used in the main proof, which guarantees that all nodes converge to the same solution when $\eta \to 0$ as $T \to \infty$.
\begin{lemma}\label{lem:bound_general}
    Given a gossip matrix as in \cref{def:gossip-matrix} (possibly extended with delays as in \cref{sec:incorporating-delays}), a large enough $\tau_g$ such that $\norm{W^{\tau_g} - W^\infty}_2^2 \leq C\rho^{\tau_g} \coloneqq 1-c < 1$, and under \cref{ass:f}, with constant stepsize $\eta$, \algname iterations satisfy
    \begin{equation*}
        \sum_{n = 1}^N \normsq{\tilde{x}(t) -  x_n(t)}_2 \leq \eta^2 \frac{4 N \norm{D}_2^2 G^2}{c^2},
    \end{equation*}
    where $D$ is the diagonal correction matrix, with virtual nodes having 1 in their corresponding entries.
\end{lemma}
\begin{IEEEproof}
    We follow the proof from \cite[Lemma A.2]{decentralized_DL_with_compression}.
    Define $r_t = \expec{\normsq{X(t) - \tilde X(t)}_F}$.
    Then,
    $$r_{t+1} = \expec{\normsq{(W^{\tau_g} - W^\infty) (X(t) -\eta D \partial F(t) - \tilde X(t))}_F},$$ since $(W^{\tau_g} - W^\infty)\tilde X(t) = 0$.
    Using the contractive property of $W^{\tau_g} - W^\infty$, it follows immediately that
    \begin{align*}
        r_{t+1} & \leq (1-c) \expec{\normsq{X(t) - \tilde X(t) - \eta D \partial F(t) }_F}    \\
                & \leq (1-c) (1+\alpha) \expec{\normsq{X(t) - \tilde X(t)}_F}                 \\
                & + (1-c) (1+\alpha^{-1})\eta^2 \expec{\normsq{D \partial F(t)}_F}            \\
                & \leq \left(1-\frac{c}{2}\right) r_t + \frac{2}{c}\eta^2 \normsq{D}_2 N G^2.
    \end{align*}
    Finally, we can check by induction that $r_t \leq \frac{\eta^2 4 \normsq{D}_2 N G^2}{c^2}$ satisfies the recursion and the proof is complete.
\end{IEEEproof}

\section{Experimental results}
\begin{figure}[htbp]
    \centering
    \includegraphics[width=.9\columnwidth]{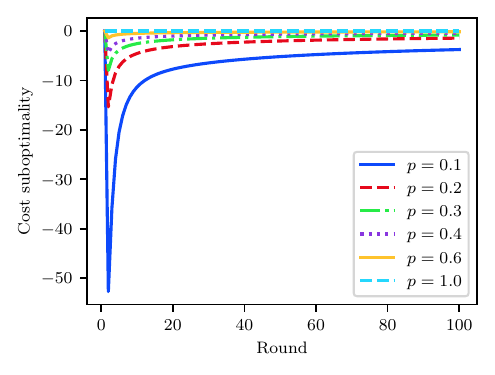}
    \caption{\algname's cost comparison for strongly connected directed random graphs with edge probability $p$. The cost is $\frac{1}{N}\sum_{n=1}^N (n-x_n(t))^2$, and the suboptimality is the difference with respect to the baseline, $p=1$, which is centralized SGD.}
    \label{fig:p-cost}
\end{figure}

\begin{figure}[htbp]
    \centering
    \includegraphics[width=.9\columnwidth]{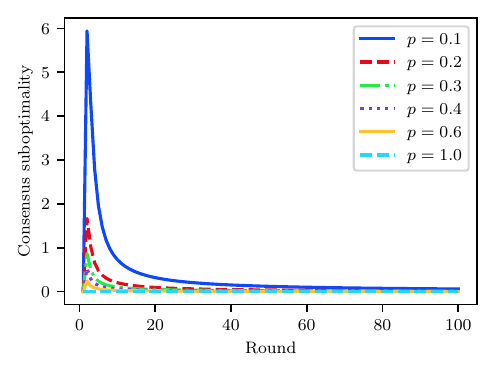}
    \caption{\algname's consensus comparison for strongly connected directed random graphs with edge probability $p$. Consensus is defined as $\frac{1}{N}\sum_{n=1}^N (\bar x(t) -x_n(t))^2$, and the suboptimality is the difference with respect to the baseline, $p=1$, which is centralized SGD.}
    \label{fig:p-consensus}
\end{figure}

We perform a set of experiments on the following problem:
Each node $n$ has an objective function $f_n(x) = (x-n)^2$.
This defines a function $f$ that has a unique global minimum at $x^*=\frac{1+2+\cdots+N}{N}$.
We generate random directed graphs with $N=100$ nodes. Every possible edge in the graph exists with a probability $p$ and non-strongly connected generated graphs are discarded. In this section, we report \algname's mean performance over $I=100$ randomly generated graphs.
Note that if $p=1$, the graph is complete and \algname is equivalent to centralized SGD, which serves as our benchmark.
Every node averages all received messages with equal weights, i.e., the inverse of the in-degree.
We use $\eta_t = \frac{0.1}{\sqrt{t}}$ and $\tau_g=1$, and report the cost and consensus suboptimalities.
The cost and consensus are defined as $\frac{1}{N}\sum_{n=1}^N (n-x_n(t))^2$ and $\frac{1}{N}\sum_{n=1}^N (\bar x(t) -x_n(t))^2$, respectively.
For both measures, the Suboptimality is defined as the difference between the value and that of the centralized SGD method, which is equivalent to the case of $p=1$.
\cref{fig:p-cost,fig:p-consensus} show how the graph topology, measured by the edge probability $p$ affects \algname's performance.
Higher graph connectivity implies faster consensus.
Also, a smaller consensus value allows solutions that differ more at each node, which make local costs lower.
In turn, the cost suboptimality has negative peaks before reaching consensus.

\begin{figure}[htbp]
    \centering
    \includegraphics[width=.9\columnwidth]{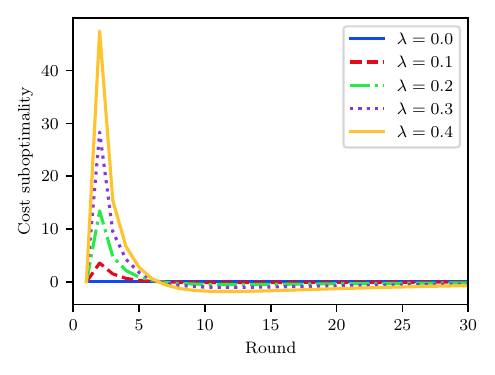}
    \caption{\algname's cost comparison for complete graphs with Poisson-distributed delays among the edges. The cost is $\frac{1}{N}\sum_{n=1}^N (n-x_n(t))^2$, averaging only over non-virtual nodes. The suboptimality is the difference with respect to the baseline, $\lambda=0$, which is centralized SGD without delays.}
    \label{fig:lambda-cost}
\end{figure}

\begin{figure}[htbp]
    \centering
    \includegraphics[width=.9\columnwidth]{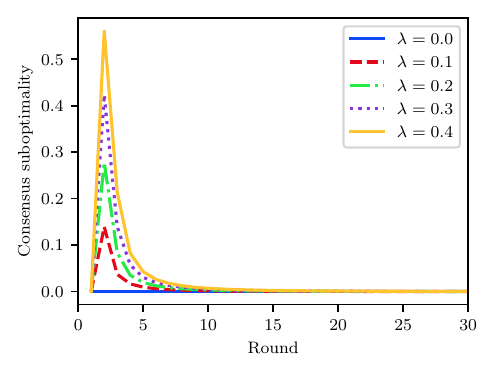}
    \caption{\algname's cost comparison for complete graphs with Poisson-distributed delays among the edges. Consensus is $\frac{1}{N}\sum_{n=1}^N (\bar x(t) -x_n(t))^2$, averaging only over non-virtual nodes. The suboptimality is the difference with respect to the baseline, $\lambda=0$, which is centralized SGD without delays.}
    \label{fig:lambda-consensus}
\end{figure}

\cref{fig:lambda-cost,fig:lambda-consensus} show the effects of delay on \algname's performance.
We consider a complete graph of $N=100$ nodes and add delays at each node according to a Poisson distribution with parameter $\lambda$.
We generate the delays for $I=100$ cases for each $\lambda$ value and report the mean performance.
Note that $\lambda=0$ implies no delays, and for a complete graph without delays \algname is the same as centralized SGD.
This will serve as the baseline for our experiments.
As expected, higher values of $\lambda$ make \algname converge slower.
\cref{fig:lambda-cost,fig:lambda-consensus} only include the first 30 iterations of the algorithm to highlight the delays' effects.
The code for both experiments is available at \cite{code_for_DT-GO}.

\section{Conclusions}
We have presented \algname, a decentralized algorithm for solving optimization problems in networks with arbitrary delays.
Nodes in the network need not have knowledge of their out-degree for \algname to converge.
We have shown that, under standard assumptions, \algname is guaranteed to converge to the problem's optimal solution at the same complexity order as the plain decentralized SGD with a fully connected graph.
We have also provided simulation results that back up our theoretical findings.

\bibliography{references}
\bibliographystyle{IEEEtran}

\end{document}